\documentclass [11pt,twoside,a4paper]{article}
\usepackage{amsfonts}
\usepackage{amsthm}
\usepackage{amsmath}
\usepackage{amstext}
\usepackage{amssymb}
\usepackage{mathrsfs}
\usepackage{amscd}
\usepackage{xypic}
\usepackage{epsf}              
\usepackage{graphicx}          
\usepackage{fancybox}          
\usepackage{color}             
\usepackage{fancyhdr}
\usepackage[hang,footnotesize]{caption2}  
\def\mathscr{\mathcal}
\newfont{\aaa}{cmb10 at 19pt}
\newfont{\bbb}{cmb10 at 11pt}

\def\v1{\vspace{1mm}}

\def\leq{\leqslant}
\def\ge{\geqslant}
\def\geq{\geqslant}

\def\dfrac{\displaystyle\frac}

\newcommand{\beq}{\begin{equation}}
	\newcommand{\eeq}{\end{equation}}
\newcommand{\bey}{\begin{eqnarray}}
	\newcommand{\eey}{\end{eqnarray}}
\newcommand{\beyy}{\begin{eqnarray*}}
	\newcommand{\eeyy}{\end{eqnarray*}}

\renewcommand{\topmargin}{0.1cm}
\renewcommand{\baselinestretch}{0.95}
\renewcommand{\headsep}{0.4cm}
\makeatletter
\def\@evenhead{
	\vbox{\hbox to \textwidth {}{\hspace{0mm}{\footnotesize
				\thepage}}{\hspace{8cm} {\footnotesize {Guoxia Feng and Manli Song}}}
		\protect\vspace{1truemm}\relax \hrule depth0pt
		height0.15truemm width\textwidth}}
\def\@evenfoot{}
\def\@oddhead{\vbox{\hbox to \textwidth
		{{\hspace{0cm}{\footnotesize Restriction theorems and Strichartz inequalities}\hfill{\footnotesize
					\thepage}}\hspace{0mm}}{} \protect\vspace{1truemm}\relax\hrule
		depth0pt height0.15truemm width\textwidth}}
\def\@oddfoot{}
\makeatother

\begin{document}
	
	
	
	\setcounter{page}{1}
	\qquad\\[8mm]
	
	\noindent{\aaa{Restriction theorems and Strichartz inequalities for the Laguerre operator involving orthonormal functions}}\\[1mm]
	
	\noindent{\bbb Guoxia Feng and Manli Song}\\[-1mm]
	
	\noindent\footnotesize{School of Mathematics and Statistics, Northwestern Polytechnical University, Xi'an, Shaanxi 710129, China}\\[6mm]
	
	\normalsize\noindent{\bbb Abstract}\quad In this paper, we prove restriction theorems for the Fourier-Laguerre transform and establish Strichartz estimates for the Schr\"{o}dinger propagator $e^{-itL_\alpha}$ for the Laguerre operator $L_\alpha=-\Delta-\sum_{j=1}^{n}(\dfrac{2\alpha_j+1}{x_j}\dfrac{\partial}{\partial x_j})+\dfrac{|x|^2}{4}$, $\alpha=(\alpha_1,\alpha_2,\cdots,\alpha_n)\in{(-\frac{1}{2},\infty)^n}$ on $\mathbb{R}_+^n$ involving systems of orthonormal functions. 
	\vspace{0.3cm}
	
	\footnotetext{
		Corresponding author: Manli Song, E-mail:
		mlsong@nwpu.edu.cn}
	
	\noindent{\bbb Keywords}\quad Laguerre operators, Strichartz inequalities, Restricction theorem.\\
	{\bbb MSC}\quad {35Q41, 47B10, 35P10, 35B65}\\[0.4cm]
	
	\newtheorem{theorem}{Theorem}[section]
	\newtheorem{preliminaries}{Preliminaries}[section]
	\newtheorem{definition}{Definition}[section]
	\newtheorem{main result}{Main Result}[section]
	\newtheorem{lemma}{Lemma}[section]
		\newtheorem{problem}{Problem}
	\newtheorem{proposition}{Proposition}[section]
	\newtheorem{corollary}{Corollary}[section]
	\newtheorem{remark}{Remark}[section]
	{\centering\section{Introduction}}   
	\subsection{Spectral of the Laguerre operator}
For any $\alpha>-1$ and $k\in\mathbb{N}$, the one-dimensional Laguerre polynomial $L_k^\alpha(x)$ of type  $\alpha$ on $\mathbb{R}_+=(0,\infty)$ is defined by 
\begin{equation*}
e^{-x}x^\alpha L_k^\alpha (x)=\frac{1}{k!}\frac{d^k}{dx^k}\left(e^{-x}x^{k+\alpha}\right), \,\forall x>0.
\end{equation*}
Every $L_k^\alpha(x)$ is a polynomial of degree $k$ which is explicitly given by
\begin{equation*}
L_k^\alpha(x)=\sum_{j=0}^{k}\dfrac{\Gamma(k+\alpha+1)}{\Gamma(k-j+1)\Gamma(j+\alpha+1)}\dfrac{(-x)^j}{\Gamma(j+1)}.
\end{equation*}
The Laguerre polynomials satisfy the following generating function identity: for $|w|<1$
\begin{equation*}
\sum_{k=0}^{\infty}L_k^\alpha(x)w^k=(1-w)^{-\alpha-1}e^{-\frac{xw}{1-w}}.
\end{equation*}
Define the Laguerre function $\psi_k^\alpha(x)$ on $\mathbb{R_+}$ by
\begin{equation*}
\psi_k^\alpha(x)={\left(\frac{2^{-\alpha}\Gamma(k+1)}{\Gamma(k+\alpha+1)}\right)}^{\frac{1}{2}}L_k^\alpha(\frac{x^2}{2})e^{-\frac{x^4}{4}}.
\end{equation*}
The Laguerre operator on $\mathbb{R}_+$ is given by
\begin{equation*}
L_\alpha=-\frac{d^2}{dx^2}-\frac{2\alpha+1}{x}\frac{d}{dx}+\frac{x^2}{4}.
\end{equation*}
Then the one-dimensional Laguerre functions $\{{\psi}_k^\alpha\}_{k=0}^\infty$ form a complete orthonormal system in $L^2(\mathbb{R}_+,dw_\alpha)$ where $dw_\alpha(x)=x^{2\alpha+1}dx$, and ${\psi}_k^\alpha$ is an eigenfunction of the Laguerre operator $L_\alpha$ with the eigenvalue $(2k+\alpha+1)$, i.e.,
\begin{equation*}
L_\alpha{\psi}_k^\alpha=(2k+\alpha+1){\psi}_k^\alpha.
\end{equation*}
Therefore, for any $f\in{L^2(\mathbb{R}_+, dw_\alpha)}$, it has the Laguerre expansion
\begin{equation*}
f=\sum_{k=0}^{\infty}\langle f,\psi_k^\alpha\rangle_\alpha\psi_k^\alpha,
\end{equation*}
where $\langle\cdot,\cdot\rangle_\alpha$ is an inner product inherited from $L^2(\mathbb{R}_+,dw_\alpha)$. It can be easily checked the Laguerre operator $L_\alpha$ is  positive and essentially self-adjoint.

For any multi-index $\mu=(\mu_1,\mu_2,\dots,\mu_n)\in\mathbb{N}^n$ and $\alpha=(\alpha_1,\alpha_2,\dots,\alpha_n)\in{(-1,\infty)^n}$, the $n$-dimensional Lagurre function $\psi_\mu^\alpha(x)$ on $\mathbb{R}_+^n$ is given by the tensor product of one-dimensional Laguerre functions
\begin{equation*}
\psi_\mu^\alpha(x)=\prod_{j=1}^{n}\psi_{\mu_j}^{\alpha_j}(x_j),\, \forall x=(x_1,x_2,\cdots,x_n)\in{\mathbb{R}_+^n},
\end{equation*}
and the $n$-dimensional Laguerre operator $L_\alpha$ is defined as the sum of one-dimensional Laguerre operator $L_{\alpha_j}$
\begin{equation*}\label{Ldef}
L_{\alpha}=\sum_{j=1}^{n}{L_{\alpha_j}}=-\Delta-\sum_{j=1}^{n}(\dfrac{2\alpha_j+1}{x_j}\dfrac{\partial}{\partial x_j})+\dfrac{|x|^2}{4}.
\end{equation*}
$L_{\alpha}$ is also positive and essentially self-adjoint and $\{\psi_\mu^\alpha\}_{\mu\in\mathbb{N}^n}$ also form a complete orthonormal system in $L^2(\mathbb{R}_+^n,dw_\alpha)$, where $dw_\alpha(x)=x_1^{2\alpha_1+1}x_2^{2\alpha_2+1}\dots x_n^{2\alpha_n+1}dx_1dx_2\dots dx_n=x^{2\alpha+1}dx$. We also have that $\psi_\mu^\alpha$ is the eigenfunction of $L_\alpha$ with the eigenvalue $2|\mu|+\sum_{j=1}^{n}\alpha_j +n$, i,e., 
 \begin{equation*}
L_\alpha\psi_\mu^\alpha=(2|\mu|+\sum_{j=1}^{n}\alpha_j+n)\psi_\mu^\alpha, \text{ where }|\mu|=\sum_{j=1}^{n}\mu_j.
\end{equation*}
Thus, for any $f\in{L^2(\mathbb{R}_+^n,dw_\alpha)}$, it has the Laguerre expansion 
\begin{equation*}
f=\sum_{\mu\in\mathbb{N}}\langle f,\psi_\mu^\alpha\rangle_\alpha\psi_\mu^\alpha:=\sum_{j=1}^{n}P_kf,
\end{equation*}
where $P_k$ denotes the orthogonal projection operator corresponding to the eigenvalue  $2k+\sum_{j=1}^{n}\alpha_j+n$, i.e.,
\begin{equation*}
P_kf=\sum_{|\mu|=k}\langle f,\psi_\mu^\alpha\rangle_\alpha\psi_\mu^\alpha.
\end{equation*}
Then the spectral decomposition of $L_\alpha$ is explicitly given by
\begin{equation*}
L_\alpha f=\sum_{k=0}^{\infty}\left(2k+\sum_{j=1}^{n}\alpha_j+n\right)P_kf.
\end{equation*}
For $f\in{L^2(\mathbb{R}_+^n,dw_\alpha)}$, $e^{-itL_\alpha}f$ is given by
\begin{equation*}
e^{-itL_\alpha}f(x)=\sum_{k=0}^{\infty}e^{-it\left(2k+\sum\limits_{j=1}^{n}\alpha_j+n\right)}P_kf=\sum_{k=0}^{\infty}e^{-it\left(2k+\sum\limits_{j=1}^{n}\alpha_j+n\right)}\sum_{|\mu|=k}\langle f,\psi_\mu^\alpha\rangle_\alpha\psi_\mu^\alpha.
\end{equation*}
It is clear that $e^{-itL_\alpha}$ is a unitary operator on $L^2(\mathbb{R}_+^n,dw_\alpha)$ with the adjoint operator $e^{itL_\alpha}$.
\begin{remark}\label{mark}
	$e^{-itL_\alpha}f$ is periodic in $t$ if and only if \ $\sum_{j=1}^{n}\alpha_j$ is rational while $|e^{-itL_\alpha}f|$ is always periodic in $t$ with period of $\pi$.
\end{remark}
\begin{lemma}[see \cite{Sohani}] Let $z=r+it$, $r>0$ and $\alpha\in{(-\frac{1}{2},\infty)^n}$. Then $e^{-zL_\alpha}$ is an integral operator on $L^2(\mathbb{R}_+^n,dw_\alpha)$ and
	\begin{align*}
	e^{-zL_\alpha}f(x)&=\int_{\mathbb{R}_+^n}f(y)K_z^\alpha(x,y)dw_\alpha(y),\\
	K_z^\alpha(x,y)&=\sum_{k=0}^{n}e^{-z\left(2k+\sum_{j=1}^{n}\alpha_j+n\right)}\sum_{|\mu|=k}\psi_\mu^\alpha(x)\psi_\mu^\alpha(y).
	\end{align*}
Moreover, for any $r\in{(0,1]}$, we have a uniform estimate for the kernel 
	\begin{equation}\label{kernel}
      |K_z^\alpha(x,y)|\leq\dfrac{C}{|sint|^{n+\sum\limits_{j=1}^{n}\alpha_j}},
	\end{equation}
where $C>0$ depends only on $n$ and $\alpha$.
\end{lemma}
\subsection{The restriction problem  for the Laguerre operator}
The restriction theorem for the Fourier transform plays an important role in harmonic analysis as well as in the theory of partial differential equations. Given a surface $S$ in $\mathbb{R}^n$, one asks for which exponent $1\leq p\leq 2$ the Fourier transform of an $L^p$-function is square integrable on $S$, which is a classical model case of restriction problem often considered in the literature. There are two types of surfaces for which this problem has been completely settled. For smooth compact surfaces with non-zero Gauss curvature,  see the work by Stein-Tomas \cite{Stein1984, Tomas}. For quadratic surfaces, Strichartz \cite{Stri} split them into three categories (parabloid-like, cone-like or sphere-like) and gave a complete answer. In this section, we introduce the restriction problem for the Laguerre operator.

For any $f\in L^1(\mathbb{R}_+^n,dw_\alpha)$ and $\alpha\in{(-\frac{1}{2},\infty)^n}$, the Fourier-Laguerre transform of $f$ is defined by
\begin{equation*}
\hat{f}(\mu,\alpha)=\int_{\mathbb{R}_+^n}f(x)\psi_\mu^\alpha(x)dw_\alpha(x),\ ,\forall \mu\in\mathbb{N}^n.
\end{equation*}
If $f\in{L^2(\mathbb{R}_+^n,dw_\alpha)}$, then $\{\hat{f}(\mu,\alpha)\}_{\mu\in\mathbb{N}^n}\in{\ell^2(\mathbb{N}^n)}$ and satisfies the Plancherel formula
\begin{equation}\label{Plancherel1}
\|f\|_{L^2(\mathbb{R}_+^n,dw_\alpha)}=\|\{\hat{f}(\mu,\alpha)\}\|_{\ell^2(\mathbb{N}^n)}.
\end{equation}
The inverse Laguerre transform is given by 
\begin{equation*}
f(x)=\sum_{\mu\in{\mathbb{N}^n}}\hat{f}(\mu,\alpha)\psi_\mu^\alpha(x), \,\forall f\in\mathcal{S}(\mathbb{R}_+^n,dw_\alpha).
\end{equation*}

It was discovered by Strichartz \cite{Stri} that the restriction theorems for some quadradic surfaces are linked to the space-time decay estimates (also called Strichartz estimates) for some envolution equations. Consider the Schr\"{o}dinger equation associated with the Laguerre operator
\begin{equation}\label{initial}
\left\{
\begin{array}{rl}
i\partial_tu(x,t)&=L_\alpha u(x,t),\, x\in{\mathbb{R}_+^n},t\in{\mathbb{R}}\\
u(x,0)&=f(x).
\end{array}
\right.
\end{equation}
We know that if $f\in L^2(\mathbb{R}_+^n,dw_\alpha)$, the solution of this initial value problem (\ref{initial}) is given by $u(x,t)=e^{-itL_\alpha}f(x)$ and $|u(x,t)|$ is periodic in $t$ with period $\pi$.

By the idea of Strichartz \cite{Stri}, the space-time decay estimate for the solution of \eqref{initial} is also reduced to the restriction theorem on $\mathbb{R}_+^n\times\mathbb{T}$ where $\mathbb{T}=[-\pi,\pi]$. So we need to introduce the Fourier-Laguerre transform on $\mathbb{R}_+^n\times\mathbb{T}$. For any $F\in{L^1(\mathbb{T}, L^1(\mathbb{R}_+^n,dw_\alpha))}$, the Fourier-Laguerre transform is given by
\begin{equation*}
\hat{F}(\mu,\nu,\alpha)=\dfrac{1}{2\pi}\int_{\mathbb{R}_+^{n}}\int_\mathbb{T}F(x,t)\psi_\mu^\alpha(x)e^{it\nu}dtdw_\alpha (x),\,\forall \mu\in{\mathbb{N}^n}, \nu\in{\mathbb{Z}}.
\end{equation*}
If $F\in L^2(\mathbb{T}, L^2(\mathbb{R}_+^n,dw_\alpha))$, then $\{\hat{F}(\mu,\nu,\alpha)\}_{ (\mu, \nu)\in{\mathbb{N}^n}\times\mathbb{Z}}\in\ell^2(\mathbb{N}^n\times\mathbb{Z})$ and the Plancherel formula is of the form
\begin{equation}\label{Plancherel2}
\|F\|_{L^2(\mathbb{T}, L^2(\mathbb{R}_+^n,dw_\alpha))}=\sqrt{2\pi}\|\{\hat{F}(\mu,\nu,\alpha)\}\|_{\ell^2(\mathbb{N}^n\times\mathbb{Z})}.
\end{equation}
The inverse Fourier-Laguerre transform is given by
\begin{equation*}
F(x,t)=\sum_{(\mu, \nu)\in{\mathbb{N}^n}\times\mathbb{Z}}\hat{F}(\mu,\nu,\alpha)\psi_{\mu}^{\alpha}(x)e^{-it\nu}.
\end{equation*}

For any discrete surface $S$ in $\mathbb{N}^n\times\mathbb{Z}$ with respect to counting measure, define the restriction operator $\mathcal{R}_SF:=\{\hat{F}(\mu,\nu,\alpha)\}_{(\mu,\nu)\in S}$ and its dual operator (also called the Fourier-Laguerre extension operator) as 
\begin{equation*}
\mathcal{E}_S(\{\hat{F}(\mu,\nu,\alpha)\})(x,t):=\sum_{(\mu,\nu)\in{S}}\hat{F}(\mu,\nu,\alpha)\psi_{\mu}^\alpha(x)e^{-it\nu}.
\end{equation*}
We consider the following restriction problem:
\begin{problem}\label{Pro}
	For which exponents of $p, q$, $1\leq p,q\leq2$, is it true that the sequence of the Fourier-Laguerre transforms of $F\in L^q(\mathbb{T}, L^p(\mathbb{R}_+^n,dw_\alpha))$ belongs to $\ell^2(S)?$ i.e.
	\begin{align*}
	\|\mathcal{R}_SF\|_{\ell^2(S)}\leq C\|F\| _{L^q(\mathbb{T}, L^p(\mathbb{R}_+^n,dw_\alpha))} ?
	\end{align*}
\end{problem}
The duality argument shows that it is completely equivalent to the boundedness of the extension operator $\mathcal{E}_S$ from $\ell^2(S)$ to $L^{q^\prime}(\mathbb{T}, L^{p^\prime}(\mathbb{R}_+^n,dw_\alpha))$ where $\frac{1}{p}+\frac{1}{p^\prime}=1$ and $\frac{1}{q}+\frac{1}{q^\prime}=1$. 

\begin{problem}\label{Pro}
	For which exponents of $p, q$, $1\leq p, q\leq2$, is it true that the extension Fourier-Laguerre operator $\mathcal{E}_S$ is bounded from $\ell^2(S)$ to $L^{q^\prime}(\mathbb{T}, L^{p^\prime}(\mathbb{R}_+^n,dw_\alpha))$ ? i.e.
	\begin{align*}
	\|\mathcal{E}_S(\{\hat{F}(\mu,\nu,\alpha)\})\|_{L^{q^\prime}(\mathbb{T}, L^{p^\prime}(\mathbb{R}_+^n,dw_\alpha))}\leq C\|\{\hat{F}(\mu,\nu,\alpha)\}\| _{\ell^2(S)} ?
	\end{align*}
\end{problem}

Note the fact that $\mathcal{E}_S$ is bounded from $\ell^2(S)$ to $L^{q^\prime}(\mathbb{T}, L^{p^\prime}(\mathbb{R}_+^n,dw_\alpha))$ if and only if $\mathcal{T}_S:=\mathcal{E}_S(\mathcal{E}_S)^*$ is bounded from $L^q(\mathbb{T}, L^p(\mathbb{R}_+^n,dw_\alpha))$ to $L^{q^\prime}(\mathbb{T}, L^{p^\prime}(\mathbb{R}_+^n,dw_\alpha))$. 
\begin{problem}
	For which exponents of $p,q$, $1\leq p,q\leq2$, the operator $\mathcal{T}_S$ is bounded from $L^q(\mathbb{T}, L^p(\mathbb{R}_+^n,dw_\alpha))$ to $L^{q^\prime}(\mathbb{T}, L^{p^\prime}(\mathbb{R}_+^n,dw_\alpha))$? i.e.
	\begin{align*}
	\|\mathcal{T}_SF\|_{L^{q^\prime}(\mathbb{T}, L^{p^\prime}(\mathbb{R}_+^n,dw_\alpha))}\leq C\|F\| _{L^q(\mathbb{T}, L^p(\mathbb{R}_+^n,dw_\alpha))} ?
	\end{align*}
\end{problem}
\begin{remark}
Throughout this paper, the mixed norm $\|\cdot\|_{L^q(\mathbb{T}, L^p(\mathbb{R}_+^n,dw_\alpha))}$ is defined by
\begin{equation*}
\|f\|_{L^q(\mathbb{T}, L^p(\mathbb{R}_+^n,dw_\alpha))}=\left(\int_{-\pi}^\pi \left(\int_{\mathbb{R}_+^n}|f(x,t)|^pdw_\alpha(x)\right)^\frac{q}{p}dt\right)^\frac{1}{q}.
\end{equation*}
\end{remark}


\subsection{Relation between estriction estimates and Strichartz inequalities for the Schr\"{o}dinger equations associated with the Laguerre operator}\label{associated}
Restriction estimates have several applications, from spectral theory to number theory. An important application of the restriction theorem for the type of sphere-like surface proved in \cite{Stri} is known as the Strichartz inequality for the solution to the Schr\"{o}dinger equation, which is a widely used tool to study the nonlinear  Schr\"{o}dinger equation.  In this section, we explain the relation between the restriction theorem on certain dicrete surface and the Strichartz estimate for the solution of the Schr\"{o}dinger equation \eqref{initial} associated with the Laguerre operator. 

Choose $S$ to be the discrete surface $S=\{(\mu,\nu)\in{\mathbb{N}^n\times\mathbb{Z}}:\nu=2|\mu|+n\}$. 
For some measurable function $f$ on $\mathbb{R}_+^n$, let
\[ \hat{F}(\mu,\nu,\alpha)=\left\{
\begin{array}{rl}
\hat{f}(\mu,\alpha), & \ \ \ \text{if }   \nu=2|\mu|+n\\
0,\ \ \ \  & \ \ \ \text{otherwise} .
\end{array}
\right. \]
Then for any $f\in L^2(\mathbb{R}_+^n,dw_\alpha)$, by the Plancherel formula \eqref{Plancherel1} and \eqref{Plancherel2}, we have
\begin{equation}\label{F-f}
\begin{aligned}
\|F\|_{(L^2(\mathbb{T}), L^2(\mathbb{R}_+^n,dw_\alpha))}&=\sqrt{2\pi}\|\{\hat{F}(\mu,\nu,\alpha)\}\|_{\ell^2(S)}\\
&=\sqrt{2\pi}\|\{\hat{f}(\mu,\alpha)\}\|_{\ell^2(\mathbb{N}^n)}=\sqrt{2\pi}\|f\|_{L^2(\mathbb{R}_+^n,dw_\alpha)}.
\end{aligned}
\end{equation}
Hence, $F\in L^2(\mathbb{T},L^2(\mathbb{R}_+^n,dw_\alpha))$ 
and 
\begin{equation}\label{relation}
\begin{aligned}
\mathcal{E}_S(\{\hat{F}(\mu,\nu,\alpha)\})(x,t)  & =\sum_{(\mu,\nu)\in{S}}\hat{F}(\mu,\nu,\alpha)\psi_\mu^\alpha(x)e^{-it\nu} \\
& =\sum_{\mu\in\mathbb{N}^n}\hat{f}(\mu,\alpha)\psi_\mu^\alpha(x)e^{-it(2|\mu|+n)} \\
& =\sum_{\mu\in\mathbb{N}^n}\langle f,\psi_\mu^\alpha\rangle_\alpha\psi_\mu^\alpha(x)e^{-it(2|\mu|+n)}\\
&=e^{it\sum_{j=1}^{n}\alpha_j}\sum_{\mu\in\mathbb{N}^n}e^{-it\left(2|\mu|+\sum_{j=1}^n\alpha_j+n\right)}\langle f,\psi_\mu^\alpha\rangle_\alpha\psi_\mu^\alpha(x)\\
&=e^{it\sum_{j=1}^{n}\alpha_j}\cdot e^{-itL_\alpha}f(x).
\end{aligned}
\end{equation}
If $\mathcal{E}_S$ is bounded from $\ell^2(S)$ to $L^{q^\prime}(\mathbb{T}, L^{p^\prime}(\mathbb{R}_+^n,dw_\alpha))$ for some $1\leq p,q\leq2$, the Strichartz inequality follows from \eqref{F-f} and \eqref{relation} that
\begin{equation*}
\begin{aligned}
\|e^{-itL_\alpha}f\|_{L^{q^\prime}(\mathbb{T}, L^{p^\prime}(\mathbb{R}_+^n,dw_\alpha))}&=\|\mathcal{E}_S(\{\hat{F}(\mu,\nu,\alpha)\})\|_{L^{q^\prime}(\mathbb{T}, L^{p^\prime}(\mathbb{R}_+^n,dw_\alpha))}\\
&\leq C\|\{\hat{F}(\mu,\nu,\alpha)\}\|_{\ell^{2}(S)}\\
&= C\sqrt{2\pi}\|f\|_{L^2(\mathbb{R}_+^n,dw_\alpha)}.
\end{aligned}
\end{equation*}
Indeed, the Strichartz inequality for the solution of the schr\"{o}dinger equation \eqref{initial} associated with the Laguerre operator holds if only if the restriction theorem holds on the specific surface $S=\{(\mu,\nu)\in{\mathbb{N}^n\times\mathbb{Z}}: \nu=2|\mu|+n\}$. 

Sohani \cite{Sohani} established Strichartz estimate for the Schr\"{o}dinger progagator $e^{-itL_\alpha}$.
\begin{lemma}\label{Original-Strichartz}
Let $n\geq1$ and $\alpha\in(-\frac{1}{2},\infty)^n$. If $1\leq p\leq \infty$, $1<q\leq \infty$ and
\begin{equation*}
\frac{1}{q}+\frac{n+\sum_{j=1}^n\alpha_j}{p}=n+\sum_{j=1}^n\alpha_j,
\end{equation*}
we have
\begin{equation*}
\|e^{-itL_\alpha}f\|_{L^{2q}(\mathbb{T}, L^{2p}(\mathbb{R}_+^n,dw_\alpha))}\leq C\|f\|_{L^2(\mathbb{R}_+^n,dw_\alpha)}.
\end{equation*}
\end{lemma}

Equivalently, we obtain the following restriction theorem on the specific surface $S=\{(\mu,\nu)\in{\mathbb{N}^n\times\mathbb{Z}}: \nu=2|\mu|+n\}$. 
\begin{theorem}(Restriction theorem for a single function)\label{Restriction-single}
Let $S=\{(\mu,\nu)\in{\mathbb{N}^n\times\mathbb{Z}}: \nu=2|\mu|+n\}$. Under the same hypotheses as in Lemma \ref{Original-Strichartz}, we have
\begin{equation*}
\|\mathcal{E}_S(\{\hat{F}(\mu,\nu,\alpha)\})\|_{L^{2q}(\mathbb{T}, L^{2p}(\mathbb{R}_+^n,dw_\alpha))}\leq C\|\{\hat{F}(\mu,\nu,\alpha)\}\|_{\ell^2(S)}.
\end{equation*}
\end{theorem}

\subsection{Generalizations to the system of orthonormal functions}
Generalizing functional inequlities involving a single function to the system of orthonormal functions is not a new topic. It is strongly motivated by the study many-body systems in quantum mechnics, where a system of $N$ independent fermions in $\mathbb{R}^n$ is decribed by a collection of $N$ orthonormal functions in $L^2(\mathbb{R}^n)$. For this reason, it is important to obtain functional inequalities involving a large number of orthonormal functions. 

The first example is the Lieb-Thirring inequality \cite{Lieb-Thirring} generating the Gagliardo-Nirenberg-Sobolev inequality, which is a decisive tool for proving stability of matter. In 2013, Frank, Lewin, Lieb and Seiringer \cite{FLLS} generalized the Strichartz inequality to a system of orthonormal functions. Lewin and Sabin \cite{Lewin-Sabin1, Lewin-Sabin2} applied the generalized Strichartz estimates to study the nonlinear evolution of quatum system with infinite many particles. 

Later, Frank and Sabin \cite{FS} obtained a duality principle in terms of Schatten class and generalized the theorems of Stein-Tomas and Strichartz about surfaces restrictions of Fourier transforms to systems of orthonormal functions. The advantage of the duality principle is that it allows to deduce automatically the restriction bounds for orthonormal systems from Schatten bounds for a single function if it could be obtained by a certain method based on complex interpolation. Now we explain how to get the Schatten bounds. Stein \cite{Stein1984} and Strichartz \cite{Stri} proved the boundedness of $\mathcal{T}_S$ by introducing an analytic family of operators $\{T_z\}$ in the sense of Stein defined in a strip $a\leq \Re(z)\leq b$ in the complex plane such that $\mathcal{T}_S=T_c$ for some $c\in [a,b]$.  They proved that $T_z$ is $L^2-L^2$ bounded on the line $\Re (z)=b$ and $L^1-L^\infty$ bounded of on the line $\Re (z)=a$. Using Stein's complex interpolation theorem \cite{Stein}, they deduced the $L^p-L^{p^\prime}$ boundedness of $\mathcal{T}_S$ for some $p\in [1,2]$. By H\"{o}lder's inequality, $\mathcal{T}_S$ is $L^p-L^{p^\prime}$ bounded if and if the operator $W_1\mathcal{T}_SW_2$ is $L^2-L^2$ bounded for any $W_1,W_2\in L^\frac{2p}{2-p}(\mathbb{R}^n)$. Indeed, Frank and Sabin \cite{FS} showed a stronger Schatten bound for $W_1\mathcal{T}_SW_2$ which was a more general result than $L^2-L^2$ boundedness.

Motivated by Strichartz \cite{Stri} and Frank-Sabin \cite{FS}, Mandal and Swain \cite{MS1, MS2} recently proved the restriction theorems and obtained the Strichartz estimates for systems of orthonormal functions for the Hermite and special Hermite operator.

To the best of our knowledge, the study on the restriction theorem with respect to the Fourier-Laguerre transform has not been considered in the literature so far. The aim of this paper is to establish the restriction theorem and Strichartz estimate for systems of orthonormal functions with respect to the Laguerre operator. We consider the restriction estimate of the form:

\begin{problem}\label{Pro}
For any orthonormal system $\left\{\{\hat{F_j}(\mu,\nu,\alpha)\}_{(\mu,\nu)\in\mathbb{N}^n\times\mathbb{Z}}\right\}_{j\in{J}}$ in $\ell^2(S)$ and any sequence $\{n_j\}_{j\in{J}}\subseteq\mathbb{C}$, we have 
\begin{equation*}
\bigg\|\sum_{j\in J}n_j|\mathcal{E}_S\{\hat{F_j}(\mu,\nu,\alpha)\}|^2\bigg\|_{L^{\frac{q^\prime}{2}}(\mathbb{T}, L^{\frac{p^\prime}{2}}(\mathbb{R}_+^n,dw_\alpha))}\leq C\bigg(\sum_{j\in J}|n_j|^{\lambda^\prime}\bigg)^{\frac{1}{\lambda^\prime}},
\end{equation*}		
for some $1\leq p, q\leq2$ and $\lambda\geq1$, with $C>0$ independent of $\left\{\{\hat{F_j}(\mu,\nu,\alpha)\}_{(\mu,\nu)\in\mathbb{N}^n\times\mathbb{Z}}\right\}_{j\in{J}}$ and $\{n_j\}_{j\in{J}}$, 
\end{problem}
If $S$ is the discrete surface $S=\{(\mu,\nu)\in{\mathbb{N}^n\times\mathbb{Z}}:\nu=2|\mu|+n\}$, by  \eqref{relation}, we deduce the corresponding generalized Strichartz estimate 
\begin{equation*}
\bigg\|\sum_{j\in J}n_j|e^{-itL_\alpha}f_j|^2\bigg\|_{L^{\frac{q^\prime}{2}}(\mathbb{T}, L^{\frac{p^\prime}{2}}(\mathbb{R}_+^n,dw_\alpha))}\leq C\bigg(\sum_{j\in J}|n_j|^{\lambda^\prime}\bigg)^{\frac{1}{\lambda^\prime}},
\end{equation*}		
for any orthonormal system $\left\{f_j\right\}_{j\in{J}}$ in $L^2(\mathbb{R}_+^n,dw_\alpha)$ and any sequence $\{n_j\}_{j\in{J}}\subseteq\mathbb{C}$.
\noindent\\[4mm]

The construction of this paper is as follows: In section $2$, we introduce the defintion of Schatten space, a complex interpolation method and the duality principle in terms of Schatten class. In section $3$, we eastablish the Schatten boundedness for the Fourier-Laguerre extension operator for some $\lambda_0>1$ by using the complex interpolation. In section $4$, we prove the restriction theorems and Strichartz inequalities for systems of orthonormal functions associated with the Laguerre operator. As an application, we obtain the global well-posedness for the nonlinear Laguerre-Hartree equation in Schatten space.

	{\centering{\section{Prilimaries}}}
\subsection{Schatten class}
We recall the defintion of Schatten space (see \cite{Simon}). Let $\mathcal{H}$ be a complex and separable Hilbert space. Let $T:  \mathcal{H}\rightarrow \mathcal{H}$ be a compact operator and $T^*$ denote the adjoint of $T$. The singular values of $T$ are the non-zero eigenvalues of $|T|:=\sqrt{T^*T}$, which form an at most countable set denoted by $\{s_j(T)\}_{j\in\mathbb{N}}$. For $\lambda>0$, the Schatten space $\mathcal{G}^\lambda(\mathcal{H})$ is defined as the space of all compact operators $T$ on $\mathcal{H}$ such that 
\begin{equation*}
	\left(\sum_{j=1}^{\infty}s_j(T)^\lambda\right)^{\frac{1}{\lambda}}<\infty.
\end{equation*}
When $\lambda\geq1$, $\mathcal{G}^\lambda(\mathcal{H})$ is a Banach space endowed with the Schatten $\lambda-$norm defined by
\begin{equation*}
\|T\|_{\mathcal{G}^\lambda(\mathcal{H})}=\left(\sum_{j=1}^{\infty}s_j(T)^\lambda\right)^{\frac{1}{\lambda}},\, \forall T\in \mathcal{G}^\lambda(\mathcal{H}).
\end{equation*}
In particular, when $\lambda=1$, an operator belonging to the space $\mathcal{G}^1(\mathcal{H})$ is known as Trace class operator; when $\lambda=2$, an operator belonging to the space $\mathcal{G}^2(\mathcal{H})$ is known as Hilbert-Schmidt class operator. 

\subsection{The complex interpolation method}
Let us first recall that a family of operators $\{T_z\}$ on $\mathbb{R}_+^n\times \mathbb{T}$ defined in a strip $a\leq \Re(z)\leq b$ in the complex plane is analytic in the sense of Stein if it has the following properties:

$(1)$ For each $z: a\leq \Re (z)\leq b$, $T_z$ is a linear transformation of simple functions on $\mathbb{R}_+^n\times \mathbb{T}$ (that is, functions that take on only a finite number of nonzero values on sets of finite measure on $\mathbb{R}_+^n\times \mathbb{T})$ to measurable functions on $\mathbb{R}_+^n\times \mathbb{T}$.

$(2)$ For all simple functions $F, G$ on $\mathbb{R}_+^n\times \mathbb{T}$, the map $z\rightarrow \langle G, T_zF\rangle_\alpha$ is analytic in $a<\Re (z)<b$ and continuous in $a\leq \Re (z)\leq b$.

$(3)$ Moreover, $\sup_{a\leq\lambda\leq b}\left|\langle G, T_{\lambda+is}F\rangle_\alpha\right|\leq C(s)$ for some $C(s)$ with at most a (double) exponential growth in $s$.\\

We obtain Schatten boundedness by a complex interpolation method and we refer to Proposition $1$ of \cite{FS} with approciate modifications.
\begin{lemma}\label{Schatten-bound}
	Let $\{T_z\}$ be an analytic family of operator on $ \mathbb{R}_+^n\times\mathbb{T}$ in the sense of Stein defined in the strip $-\lambda_0\leq \Re (z) \leq 0$ for some $\lambda_0>1$. Assume that we have the following bounds
	\begin{equation*}		
	\begin{aligned}
	&	\|T_{is}\|_{L^2(\mathbb{T}, L^2(\mathbb{R}_+^n,dw_\alpha))\to L^2(\mathbb{T}, L^2(\mathbb{R}_+^n,dw_\alpha))}\leq M_0e^{a|s|}, \\			
	&	\|T_{-\lambda_0+is}\|_{L^1(\mathbb{T}, L^1(\mathbb{R}_+^n,dw_\alpha))\to L^\infty(\mathbb{T}, L^\infty(\mathbb{R}_+^n,dw_\alpha))}\leq M_1e^{b|s|},
	\end{aligned}
	\end{equation*}
	for all $s\in{\mathbb{R}}$ and for some $a,b,M_0,M_1 \geq 0$. Then for all $W_1,W_2\in L^{2\lambda_0}(\mathbb{T}, L^{2\lambda_0}(\mathbb{R}_+^n,dw_\alpha))$, the operator $W_1T_{-1}W_2$ belongs to $\mathcal{G}^{2\lambda_0}\left(L^2(\mathbb{T}, L^2(\mathbb{R}_+^n,dw_\alpha))\right)$, and we have the estimate
	\begin{equation*}
	\left\|W_1T_{-1}W_2\right\|_{\mathcal{G}^{2\lambda_0}\left(L^2(\mathbb{T}, L^2(\mathbb{R}_+^n,dw_\alpha))\right)}\leq C\left\|W_1\right\|_{L^{2\lambda_0}(\mathbb{T}, L^{2\lambda_0}(\mathbb{R}_+^n,dw_\alpha))}\left\|W_2\right\|_{L^{2\lambda_0}(\mathbb{T}, L^{2\lambda_0}(\mathbb{R}_+^n,dw_\alpha))}
	\end{equation*}
	where $C=M_0^{1-\frac{1}{\lambda_0}}M_1^\frac{1}{\lambda_0}$.
\end{lemma}
\subsection{The duality principle}
In order to deduce bounds for orthonormal systems from Schatten bounds for a single function, we need the duality principle lemma in our context. We refer to Lemma 3 of \cite{FS} with approciate modifications to obtain the following result.
\begin{lemma}(Duality principle)\label{duality-principle} Let $S$ be a discrete surface on $\mathbb{N}^n\times\mathbb{Z}$ and $\lambda\geq1$. Assume that $A$ is a bounded linear operator from $\ell^2(S)$ to $L^{q^\prime}(\mathbb{T}, L^{p^\prime}(\mathbb{R}_+^n,dw_\alpha))$ for some $p,q\geq 1$. Then the following statements are equivalent.

$(1)$ There is a constant $C>0$ such that for all $W\in L^2(\mathbb{T}, L^2(\mathbb{R}_+^n,dw_\alpha))$
\begin{equation}\label{Schatten}
\big\|WAA^*\overline{W}\big\|_{\mathcal{G}^\lambda(L^2(\mathbb{T}, L^2(\mathbb{R}_+^n,dw_\alpha)))}\leq C\big\|W\big\|_{L^{\frac{2q}{2-q}}(\mathbb{T}, L^{\frac{2p}{2-p}}(\mathbb{R}_+^n,dw_\alpha))}^2,
\end{equation}
where the function $W$ is interpreted as an operator which acts by multiplication.

$(2)$ There is a constant $C'>0$ such that for any orthonormal system $\left\{\{\hat{F_j}(\mu,\nu,\alpha)\}_{(\mu,\nu)\in\mathbb{N}^n}\right\}_{j\in{J}}$ in $\ell^2(S)$ and any sequence $(n_j)_{j\in{J}}$ in $\mathbb{C}$ 
\begin{equation*}
\bigg\|\sum_{j\in J}n_j|A\{\{\hat{F_j}(\mu,\nu,\alpha)\}|^2\bigg\|_{L^{\frac{q^\prime}{2}}(\mathbb{T}, L^{\frac{p^\prime}{2}}(\mathbb{R}_+^n,dw_\alpha))}\leq C'\bigg(\sum_{j\in J}|n_j|^{\lambda^\prime}\bigg)^{\frac{1}{\lambda^\prime}}.
\end{equation*}		
\end{lemma}
\begin{remark}
Lemma \ref{Schatten-bound} and \ref{duality-principle} are also valid in the domain $\mathbb{R}_+^n\times [-\frac{\pi}{2},\frac{\pi}{2}]$.
\end{remark}

{\centering{\section{The Schatten boundedness for $\mathcal{T}_S$}}}
	In this section, we apply the complex interpolation method to establish the Schatten boundedness for some $\lambda_0>1$.
\subsection{On general discrete surface}
	Let $S$ be a discerete surface $S=\{(\mu,\nu)\in{\mathbb{N}^n\times \mathbb{Z}}: R(\mu,\nu)=0\}$, where $R(\mu,\nu)$ is a polynomial of degree one,with respect to the counting measue. 
	
	For some $\lambda_0>1$ and $-\lambda_0\leq \Re (z)\leq 0$, consider the analytic family of generalized functions $G_z(\mu,\nu)$ on $\mathbb{N}^n\times\mathbb{Z}$ defined by 
\begin{equation*}
	G_z(\mu,\nu):=\gamma(z)R(\mu,\nu)_+^z
\end{equation*}		
where
\[  
R(\mu,\nu)_+^z = \left\{
\begin{array}{rl}
R(\mu,\nu)^z,\ \ \ \ \  & \text{for } R(\mu,\nu)>0,\\	0 ,\  \ \ \ \ \ \ \ \ & \text{for } R(\mu,\nu)=0,\\
\end{array} 	\right.		\]
and $\gamma(z)$ is an approciate analytic function.

For a Schwartz class funcion $\Phi$ on $\mathbb{N}^n\times\mathbb{Z}$, we have 
\begin{equation*}
 \langle G_z,\Phi \rangle:=\varphi(z)\sum_{(\mu,\nu)\in\mathbb{N}^n\times\mathbb{Z}}R(\mu,\nu)_+^z\Phi(\mu,\nu)
\end{equation*}
and
	\begin{equation*}
		\lim_{z \rightarrow -1}
		\langle G_z,\Phi \rangle=\lim_{z \rightarrow -1}\varphi(z)\sum_{(\mu,\nu)\in\mathbb{N}^n\times\mathbb{Z}}R(\mu,\nu)_+^z\Phi(\mu,\nu)=\sum_{(\mu,\nu)\in{S}}\Phi(\mu,\nu),
	\end{equation*}	
where $\varphi(z)$ has adequate properties according to the type of discrete surface $S$ considered but in all cases it has at most exponential growth at infinity when $\Re(z)=0$ and a simple zero at $z=-1$ which ensures that $G_{-1}\equiv\delta_S$. We refer the reader to \cite{GS1} for the distribution calculus of $R(\mu,\nu)^z_+$.

Moreover, the family of operators $T_z$ on $\mathbb{R}_+^n\times\mathbb{T}$ is defined by 
\begin{equation*}
	T_zF(x,t)=\sum_{(\mu,\nu)\in\mathbb{N}^n\times\mathbb{Z}}\hat{F}(\mu,\nu,\alpha)G_z(\mu,\nu)\psi_{\mu}^\alpha(x)e^{-it\nu}.
\end{equation*}		
Then  $\{T_z\}$ is an analytic family of operator on $\mathbb{R}_+^n\times\mathbb{T}$ in the sense of Stein in the strip $-\lambda_0\leq \Re (z) \leq 0$. We have the identity $\mathcal{T}_S=T_{-1}$ and 
\begin{equation}\label{qwer}
	\begin{aligned}
		T_zF(x,t) & =\frac{1}{2\pi}\sum_{(\mu,\nu)\in\mathbb{N}^n\times\mathbb{Z}}\bigg(\int_{\mathbb{R}_+^n}\int_{-\pi}^{\pi}F(y,\tau)\psi_{\mu}^\alpha(y)e^{i\tau\nu}d\tau dw_\alpha(y)\bigg)G_z(\mu,\nu)\psi_{\mu}^\alpha(x)e^{-it\nu}	\\	&=\frac{1}{2\pi}\int_{\mathbb{R}_+^n}\int_{-\pi}^{\pi}\sum_{(\mu,\nu)\in\mathbb{N}^n\times\mathbb{Z}}\psi_{\mu}^\alpha(x)\psi_{\mu}^\alpha(y)G_z(\mu,\nu)e^{-i\nu(t-\tau)} F(y,\tau)d\tau dw_\alpha(y)\\
		&=\int_{\mathbb{R}_+^n}\Big(\mathcal{K}^\alpha_z(x,y,\cdot)\ast F(y,\cdot)\Big)(t)dw_\alpha(y)
	\end{aligned}
\end{equation}
where 
\begin{equation*}
\mathcal{K}^\alpha_z(x,y,t)=\frac{1}{2\pi}\sum_{(\mu,\nu)\in\mathbb{N}^n\times\mathbb{Z}}\psi_{\mu}^\alpha(x)\psi_{\mu}^\alpha(y)G_z(\mu,\nu)e^{-i\nu t}.
\end{equation*}
When $\Re (z)=0$, we have
\begin{equation}\label{Rez=0}
		\|T_{is}\|_{L^2(\mathbb{T}, L^2(\mathbb{R}_+^n,dw_\alpha)) \to L^2(\mathbb{T}, L^2(\mathbb{R}_+^n,dw_\alpha))}=\|G_{is}\|_{\ell^\infty(\mathbb{N}^n\times\mathbb{Z})}\leq|\varphi(is)|,
\end{equation}
and using H\"{o}lder and Young inequality, we obtain 
\begin{equation}\label{Rez=1}
\|T_z\|_{L^1(\mathbb{T}, L^1(\mathbb{R}_+^n, dw_\alpha)) \to L^\infty(\mathbb{T}, L^\infty(\mathbb{R}_+^n, dw_\alpha))}\leq \sup_{x,y\in\mathbb{R}_+^n, t\in \mathbb{T}}|\mathcal{K}^\alpha_z(x,y,t)|. 
\end{equation}

By Lemma \ref{Schatten-bound}, we immediately obtain the Schatten boundedness of the form \eqref{Schatten}.
\begin{lemma} \label{General} Let $n\geq1$ and $S\subseteq\mathbb{N}^n\times\mathbb{Z}$ be a discrete surface. Suppose that for every $x,y\in\mathbb{R}_+^n$ and $t\in\mathbb{T}$, $|\mathcal{K}^\alpha_z(x,y,t)|$ is uniformly bounded by a constant $C(s)>0$ with at most exponential growth in $s$ when $z=-\lambda_0+is$ for some $\lambda_0>1$.  Then for all $W_1,W_2\in L^{2\lambda_0}(\mathbb{T}, L^{2\lambda_0}(\mathbb{R}_+^n,dw_\alpha))$,  the operator $W_1\mathcal{T}_SW_2$=$W_1T_{-1}W_2$ belongs to $\mathcal{G}^{2\lambda_0}\left(L^2(\mathbb{T}, L^2(\mathbb{R}_+^n,dw_\alpha))\right)$ and we have the estimate
	\begin{equation*}
	\left\|W_1\mathcal{T}_SW_2\right\|_{\mathcal{G}^{2\lambda_0}\left(L^2(\mathbb{T}, L^2(\mathbb{R}_+^n,dw_\alpha))\right)}\leq C\left\|W_1\right\|_{L^{2\lambda_0}(\mathbb{T}, L^{2\lambda_0}(\mathbb{R}_+^n,dw_\alpha))}\left\|W_2\right\|_{L^{2\lambda_0}(\mathbb{T}, L^{2\lambda_0}(\mathbb{R}_+^n,dw_\alpha))},
	\end{equation*}
where $C>0$ is independent of $W_1$ and $W_2$.
\end{lemma}

\subsection{On the particular surface $S=\{(\mu,\nu)\in\mathbb{N}^n\times\mathbb{Z}:|\nu|=2|\mu|+n\}$}
From Section \ref{associated}, the Strichartz estimates of the Schr\"{o}dinger equation for the Laguerre operator is closely related to the restriction theorem on $S=\{(\mu,\nu)\in\mathbb{N}^n\times\mathbb{Z}:|\nu|=2|\mu|+n\}$. Thus, we work on this particular surface.

Setting $\varphi(z)=\dfrac{1}{\Gamma(z+1)}$ and $R(\mu,\nu)=\nu-(2|\mu|+n)$, then we have
\begin{equation*}
	G_z(\mu,\nu)=\dfrac{1}{\Gamma(z+1)} \big(\nu-(2|\mu|+n)\big)_+^z.
\end{equation*}
For a Schwartz class function $\Phi$ on $\mathbb{N}^n\times\mathbb{Z}$,
\begin{align*}
		\lim_{z \rightarrow -1}
	\langle G_z,\Phi \rangle &=\dfrac{1}{\Gamma(z+1)}\sum_{\mu,\nu\in{S}}\Phi(\mu,\nu)\big(\nu-(2|\mu|+n)\big)_+^z\\
	&= \sum_{\mu,\nu\in{S}}\Phi(\mu,\nu).	
\end{align*}
Thus, $G_{-1}=\delta_S$ and 
\begin{equation*}
	T_zF(x,t)=\int_{\mathbb{R}_+^n}\Big(\mathcal{K}_z^\alpha(x,y,\cdot)* F(y,\cdot)\Big)(t)dw_\alpha(y).
\end{equation*}
By computation, we have
\begin{equation}\label{Kernel}
	\begin{aligned}
	\mathcal{K}^\alpha_z(x,y,t) & =\frac{1}{2\pi}\sum_{(\mu,\nu)\in\mathbb{N}^n\times\mathbb{Z}}\psi_\mu^\alpha(x)\psi_\mu^\alpha(y)G_z(\mu,\nu)e^{-it\nu}\\
	& = \dfrac{1}{2\pi\Gamma(z+1)}\sum_{(\mu,\nu)\in\mathbb{N}^n\times\mathbb{Z}}\psi_\mu^\alpha(x)\psi_\mu^\alpha(y)e^{-it\nu}\big(\nu-(2|\mu|+n)\big)_+^z\\
    & = \dfrac{e^{it\sum_{j=1}^{n}\alpha_j}}{2\pi\Gamma(z+1)}\sum_{\mu\in\mathbb{N}^n}\bigg(\psi_\mu^\alpha(x)\psi_\mu^\alpha(y)e^{-it(2|\mu|+\sum_{j=1}^{n}\alpha_j+n)}\bigg)\sum_{k=0}^{n}k_+^ze^{-itk}.\\
    & = \dfrac{e^{it\sum_{j=1}^{n}\alpha_j}}{2\pi\Gamma(z+1)}K^\alpha_z(x,y)\sum_{k=0}^{n}k_+^ze^{-itk},
	\end{aligned}
\end{equation}
where the last equality comes from the spectral decomposition of $e^{-itL_\alpha}$.

In order to obtain a uniformly estimate of the kernel $\mathcal{K}^\alpha_z$, we need to calculate the series $\sum_{k=0}^{n}k_+^ze^{-itk}$.
\begin{lemma}(see \cite{MS1, MS2})\label{Intermediate}
Let	$-\lambda_0\leq \Re (z)\leq 0$ for some $\lambda_0>1$. Then the series $\sum\limits_{k=0}^{\infty}k_+^ze^{-itk}$ is the Fourier seriers of a integrable functions on $[-\pi,\pi]$ which is of class $C^\infty$ on $[-\pi,\pi] \setminus \{0 \}$. Near origin this function has the same singularity as the function whose values are $\Gamma(z+1)(it)^{-z-1}$, i.e.,
\begin{equation*}
	\sum_{k=0}^{\infty}k_+^ze^{-itk} \thicksim \Gamma(z+1)it^{-z-1}+b(t),
\end{equation*}
where $b\in{C^\infty}[-\pi,\pi].$
\end{lemma}

When $\Re (z)=0$, from \eqref{Rez=0}, we have
\begin{equation}\label{Exp1}
\begin{aligned}
	\|T_{is}\|_{L^2(\mathbb{T}, L^2(\mathbb{R}_+^n,dw_\alpha)) \to L^2(\mathbb{T}, L^2(\mathbb{R}_+^n,dw_\alpha))} \leq|\varphi(is)|
=\dfrac{1}{|\Gamma(1+is)|}\leq Ce^{\pi\frac{|s|}{2}}.\\
\end{aligned}
\end{equation}
When $z=-\lambda_0+is$, from \eqref{Rez=1}, we know that $T_z$ is bounded from $L^1(\mathbb{T}, L^1(\mathbb{R}_+^n, dw_\alpha))$ to $L^\infty(\mathbb{T}, L^\infty(\mathbb{R}_+^n, dw_\alpha))$ if $|\mathcal{K}^\alpha_z(x,y,t)\big|$ is uniformly bounded for each $x,y\in{\mathbb{R}_+^n}$ and $t\in \mathbb{T}$. Since it follows from \eqref{kernel}, \eqref{Kernel} and Lamma \ref{Intermediate} that
\begin{equation}\label{KER}
\begin{aligned}
\big|\mathcal{K}^\alpha_z(x,y,t)\big|
&=\left|\frac{K_z^\alpha(x,y)}{2\pi\Gamma(z+1)}\sum_{k=0}^{n}k_+^ze^{-itk}\right|\\
&\leq\frac{C}{|t|^{Re(z+n+\sum_{j=1}^{n}\alpha_j+1)}},\,\forall x,y\in{\mathbb{R}_+^n}, \forall t\in \mathbb{T}.
\end{aligned}
\end{equation}
Therefore, for every $x,y\in{\mathbb{R}_+^n}$ and $t\in \mathbb{T}$, $|\mathcal{K}^\alpha_z(x,y,t)|$ is uniformly bounded if $\Re (z)=-(n+\sum_{j=1}^{n}\alpha_j+1)$.

Choosing $\lambda_0=n+\sum_{j=1}^{n}\alpha_j+1>1$ and then we have
\begin{equation}\label{Exp2}
\|T_{-\lambda_0+is}\|_{L^1(\mathbb{T}, L^1(\mathbb{R}_+^n, dw_\alpha)) \to L^\infty(\mathbb{T}, L^\infty(\mathbb{R}_+^n, dw_\alpha))}\leq \sup_{x,y\in\mathbb{R}_+^n, t\in \mathbb{T}}|\mathcal{K}^\alpha_{-\lambda_0+is}(x,y,t)|\leq C. 
\end{equation}

Now the Schatten boundedness comes out from Theorem \ref{General}, \eqref{Exp1} and \eqref{Exp2}. 
\begin{theorem}(Schatten bound for the extension operator)\label{Schatten-Bound}
Let $S=\{(\mu,\nu)\in\mathbb{N}^n\times\mathbb{Z}:|\nu|=2|\mu|+n\}$ be the discrete surface on $\mathbb{N}^n\times\mathbb{Z}$ and $\lambda_0=n+\sum_{j=1}^{n}\alpha_j+1$. Then for all $W_1,W_2\in L^{2\lambda_0}(\mathbb{T}, L^{2\lambda_0}(\mathbb{R}_+^n,dw_\alpha))$, the operator $W_1\mathcal{T}_SW_2$=$W_1T_{-1}W_2$ belongs to $\mathcal{G}^{2\lambda_0}\left(L^2(\mathbb{T}, L^2(\mathbb{R}_+^n,dw_\alpha))\right)$ and we have the estimate
\begin{equation*}
\left\|W_1\mathcal{T}_SW_2\right\|_{\mathcal{G}^{2\lambda_0}\left(L^2(\mathbb{T}, L^2(\mathbb{R}_+^n,dw_\alpha))\right)}\leq C\left\|W_1\right\|_{L^{2\lambda_0}(\mathbb{T}, L^{2\lambda_0}(\mathbb{R}_+^n,dw_\alpha))}\left\|W_2\right\|_{L^{2\lambda_0}(\mathbb{T}, L^{2\lambda_0}(\mathbb{R}_+^n,dw_\alpha))}.
\end{equation*} 
\end{theorem}

{\centering{\section{Restriction theorems and Strichartz inequalities for  orthonomal functions}}}
Let $S=\{(\mu,\nu)\in\mathbb{N}^n\times\mathbb{Z}:|\nu|=2|\mu|+n\}$ be the discrete surface on $\mathbb{N}^n\times\mathbb{Z}$ and $\lambda_0=n+\sum_{j=1}^{n}\alpha_j+1$. From Theorem \ref{Schatten-Bound}, we have
\begin{equation}\label{Schatten-decompostion}
\begin{aligned}
&\left\|W_1\mathcal{E}_S\left(\mathcal{E}_S\right)^*W_2\right\|_{\mathcal{G}^{2\lambda_0}\left(L^2(\mathbb{T}, L^2(\mathbb{R}_+^n,dw_\alpha))\right)}\\
&\leq C\left\|W_1\right\|_{L^{2\lambda_0}(\mathbb{T}, L^{2\lambda_0}(\mathbb{R}_+^n,dw_\alpha))}\left\|W_2\right\|_{L^{2\lambda_0}(\mathbb{T}, L^{2\lambda_0}(\mathbb{R}_+^n,dw_\alpha))}.
\end{aligned}
\end{equation} 

Taking $A=\mathcal{E}_S$, by Theorem \ref{Restriction-single}, \eqref{Schatten-decompostion} and Lemma \ref{duality-principle}, we have the restriction theorem for the system of orthonomal functions.
\begin{theorem} (Restriction estimates for orthonormal functions-diagonal case) \label{diagonal-restriction} Let $n\geq 1$ and $S=\{(\mu,\nu)\in\mathbb{N}^n\times\mathbb{Z}:|\nu|=2|\mu|+n\}$.
For any (possible infinity) orthonormal system $\left\{\{\hat{F_j}(\mu,\nu,\alpha)\}_{(\mu,\nu)\in\mathbb{N}^n\times\mathbb{Z}}\right\}_{j\in{J}}$ in $\ell^2(S)$ and any sequence $\{n_j\}_{j\in{J}}$ in $\mathbb{C}$
\begin{equation*}
\begin{aligned}
&\bigg\|\sum_{j\in J}n_j|\mathcal{E}_S\{\{\hat{F_j}(\mu,\nu,\alpha)\}|^2\bigg\|_{L^{1+\frac{1}{n+\sum_{j=1}^{n}\alpha_j}}(\mathbb{T}, L^{1+\frac{1}{n+\sum_{j=1}^{n}\alpha_j}}(\mathbb{R}_+^n,dw_\alpha))}\\
&\leq C\bigg(\sum_{j}|n_j|^{\frac{2(n+\sum_{j=1}^{n}\alpha_j+1)}{2(n+\sum_{j=1}^{n}\alpha_j)+1}}\bigg)^{\frac{2(n+\sum_{j=1}^{n}\alpha_j)+1}{2(n+\sum_{j=1}^{n}\alpha_j+1)}},
\end{aligned}
\end{equation*}		
where $C>0$ is independent of $\left\{\{\hat{F_j}(\mu,\nu,\alpha)\}_{(\mu,\nu)\in\mathbb{N}^n\times\mathbb{Z}}\right\}_{j\in{J}}$ and $\{n_j\}_{j\in{J}}$.
\end{theorem}
By \eqref{relation}, we generalize the Strichartz inequality involving systems of orthonormal functions.
\begin{theorem}\label{diagonal}(Strichartz inequalities for orthonormal functions-diagonal case) Let $n\geq1$. For any (possible infinity) orthonormal system $\left\{f_j\right\}_{j\in{J}}$ in $L^2(\mathbb{R}_+^n,dw_\alpha)$ and any sequence $(n_j)_{j\in{J}}$ in $\mathbb{C}$, we have
	\begin{equation*}
	\begin{aligned}
	&\bigg\|\sum_{j\in J}n_j|e^{-itL_ \alpha}f_j|^2\bigg\|_{L^{1+\frac{1}{n+\sum_{j=1}^{n}\alpha_j}}(\mathbb{T}, L^{1+\frac{1}{n+\sum_{j=1}^{n}\alpha_j}}(\mathbb{R}_+^n,dw_\alpha))}\\
	&\leq C\bigg(\sum_{j}|n_j|^{\frac{2(n+\sum_{j=1}^{n}\alpha_j+1)}{2(n+\sum_{j=1}^{n}\alpha_j)+1}}\bigg)^{\frac{2(n+\sum_{j=1}^{n}\alpha_j)+1}{2(n+\sum_{j=1}^{n}\alpha_j+1)}},
	\end{aligned}
	\end{equation*}		
where $C>0$ is independent of $\{f_j\}_{j\in J}$ and $\{n_j\}_{j\in J}$.
\end{theorem}

To obtain the Strichartz inequality for the general case, we prove the following Schatten boundedness.
\begin{lemma}(General Schatten bound for the extension operator)\label{general-schatten-bound}
	Let $S=\{(\mu,\nu)\in\mathbb{N}^n\times\mathbb{Z}:|\nu|=2|\mu|+n\}$.Then for any $p,q\ge1$ satisfying
	\begin{equation*}
		\frac{1}{q}+\frac{n+\sum_{j=1}^{n}\alpha_j}{p}=\frac{1}{2},\quad 2(n+\sum_{j=1}^{n}\alpha_j)+1<p\leq2(n+\sum_{j=1}^{n}\alpha_j+1) ,
	\end{equation*}
we have 
\begin{equation*}
\left\|W_1\mathcal{T}_SW_2\right\|_{\mathcal{G}^p\left(L^2(\mathbb{T}, L^2(\mathbb{R}_+^n,dw_\alpha))\right)}\leq C\left\|W_1\right\|_{L^q(\mathbb{T}, L^p(\mathbb{R}_+^n,dw_\alpha))}\left\|W_2\right\|_{L^q(\mathbb{T}, L^p(\mathbb{R}_+^n,dw_\alpha))},
\end{equation*} 
for all $W_1,W_2\in L^q(\mathbb{T}, L^p(\mathbb{R}_+^n,dw_\alpha))$, where $C>0$ is independent of $W_1,W_2$.
\end{lemma}
{\bf Proof.} For $0<\lambda<\lambda_0=n+\sum_{j=1}^{n}\alpha_j$, the operator $T_{-\lambda+is}$ is an integral operator with kernel $\mathcal{K}^\alpha_{-\lambda+is}(x,t,t-\tau)$ defined in (\ref{Kernel}), Applying the Hardy-Littlewood-Sobolev inequality along with (\ref{Exp1}) and (\ref{KER}), we have  
\begin{equation*}
	\begin{aligned}
&\|W_1^{\lambda-is}T_{-\lambda+is}W_2^{\lambda-is}\|_{\mathcal{G}^2\left(L^2([-\frac{\pi}{2},\frac{\pi}{2}], L^2(\mathbb{R}_+^n,dw_\alpha))\right)}^2\\
		  =&\int^{\frac{\pi}{2}}_{-\frac{\pi}{2}}\int^{\frac{\pi}{2}}_{-\frac{\pi}{2}}\int_{\mathbb{R}_+^{2n}}|W_1(x,t)|^{2\lambda}|\mathcal{K}^\alpha_{-\lambda+is}(x,y,t-\tau)|^2|W_2(y,\tau)|^{2\lambda}dxdydtd\tau\\
	\leq  &C\int^{\frac{\pi}{2}}_{-\frac{\pi}{2}}\int^{\frac{\pi}{2}}_{-\frac{\pi}{2}}\int_{\mathbb{R}_+^{2n}}\frac{|W_1(x,t)|^{2\lambda}|W_2(y,\tau)|^{2\lambda}}{|t-\tau|^{2(n+\sum_{j=1}^{n}\alpha_j+1)-2\lambda}}dxdydtd\tau\\
	 \leq  &C\int^{\frac{\pi}{2}}_{-\frac{\pi}{2}}\int^{\frac{\pi}{2}}_{-\frac{\pi}{2}}\frac{\|W_1(\cdot, t)\|_{L^{2\lambda}(\mathbb{R}_+^n,dw_\alpha)}^{2\lambda}\|W_2(\cdot, \tau)\|_{L^{2\lambda}(\mathbb{R}_+^n,dw_\alpha)}^{2\lambda}}{|t-\tau|^{2(n+\sum_{j=1}^{n}\alpha_j+1)-2\lambda}}dtd\tau\\
	\leq  &C\|W_1\|_{L^{\frac{2\lambda}{\lambda-n-\sum_{j=1}^{n}\alpha_j}}\left([-\frac{\pi}{2},\frac{\pi}{2}], L^{2\lambda}(\mathbb{R}_+^n,dw_\alpha)\right)}^{2\lambda}\|W_2\|_{L^{\frac{2\lambda}{\lambda-n-\sum_{j=1}^{n}\alpha_j}}\left([-\frac{\pi}{2},\frac{\pi}{2}], L^{2\lambda}(\mathbb{R}_+^n,dw_\alpha)\right)}^{2\lambda}
	\end{aligned}
\end{equation*}               
  where required $0\leq2(n+\sum_{j=1}^{n}\alpha_j+1)-2\lambda<1$, that is $2(n+\sum_{j=1}^{n}\alpha_j)+1<2\lambda\leq2(n+\sum_{j=1}^{n}\alpha_j+1)$. By Theorem 2.9 of \cite{Simon} and the identity $\mathcal{T}_S=T_{-1}$, we have 
  \begin{equation}\label{part}
  	\begin{aligned}
  	&\left\|W_1 \mathcal{T}_SW_2\right\|_{\mathcal{G}^{2\lambda}\left(L^2([-\frac{\pi}{2},\frac{\pi}{2}], L^2(\mathbb{R}_+^n,dw_\alpha))\right)}=\left\|W_1 T_{-1}W_2\right\|_{\mathcal{G}^{2\lambda}\left(L^2([-\frac{\pi}{2},\frac{\pi}{2}], L^2(\mathbb{R}_+^n,dw_\alpha))\right)}\\  \leq &C\|W_1\|_{L^{\frac{2\lambda}{\lambda-n-\sum_{j=1}^{n}\alpha_j}}\left([-\frac{\pi}{2},\frac{\pi}{2}], L^{2\lambda}(\mathbb{R}_+^n,dw_\alpha)\right)}\|W_2\|_{L^{\frac{2\lambda}{\lambda-n-\sum_{j=1}^{n}\alpha_j}}\left([-\frac{\pi}{2},\frac{\pi}{2}], L^{2\lambda}(\mathbb{R}_+^n,dw_\alpha)\right).}
  	\end{aligned} 
  	\end{equation}
By Lemma \ref{duality-principle}, we have
\begin{equation}\label{1}
\begin{aligned}
	&\bigg\|\sum_{j\in J}n_j|\mathcal{E}_S\{\{\hat{F_j}(\mu,\nu,\alpha)\}|^2\bigg\|_{L^{\frac{\lambda}{n+\sum_{j=1}^{n}\alpha_j}}([-\frac{\pi}{2},\frac{\pi}{2}], L^{\frac{\lambda}{\lambda-1}}(\mathbb{R}_+^n,dw_\alpha))}\\
	&\leq C\bigg(\sum_{j}|n_j|^{(2\lambda)^\prime}\bigg)^{(2\lambda)^\prime},
\end{aligned}
\end{equation}		 
for any orthonormal system $\left\{\{\hat{F_j}(\mu,\nu,\alpha)\}_{(\mu,\nu)\in\mathbb{N}^n\times\mathbb{Z}}\right\}_{j\in{J}}$ in $\ell (S)$ and $\{n_j\}_{j\in{J}}$ in $\mathbb{C}$, which is equivalent to
\begin{equation}\label{2}
\begin{aligned}
&\bigg\|\sum_{j\in J}n_j|e^{-itL_ \alpha}f_j|^2\bigg\|_{L^{\frac{\lambda}{n+\sum_{j=1}^{n}\alpha_j}}([-\frac{\pi}{2},\frac{\pi}{2}], L^{\frac{\lambda}{\lambda-1}}(\mathbb{R}_+^n,dw_\alpha))}\\
	&\leq C\bigg(\sum_{j}|n_j|^{(2\lambda)^\prime}\bigg)^{(2\lambda)^\prime},
\end{aligned}
\end{equation}		
for any orthonormal system $\{f_j\}_{j\in J}$ in $L^2(\mathbb{R}_+^n, dw_\alpha)$ and $\{n_j\}_{j\in J}$ in $\mathbb{C}$.

Note that $|e^{-itL_ \alpha}f_j|$ is periodic with period $\pi$, so \eqref{1} and \eqref{2} are also valid for the domain $\mathbb{R}_+^n\times\mathbb{T}$. Using Lemma \ref{duality-principle} again, \eqref{part} holds true also for the domain $\mathbb{R}_+^n\times\mathbb{T}$. \qed
\\[11pt]

Now we are ready to establish the restriction theorem for systems of orthonormal functions.
\begin{theorem}(Restriction estimates for orthonormal functions-general case)\label{main-restriction}  If $p,q,n\geq1$ such that 
	\begin{equation*}
	1\leq p < \dfrac{2(n+\sum_{j=1}^{n}\alpha_j)+1}{2(n+\sum_{j=1}^{n}\alpha_j)-1}\ \ \ \  and \ \ \ \ \dfrac{1}{q}+\dfrac{n+\sum_{j=1}^{n}\alpha_j}{p}=n+\sum_{j=1}^{n}\alpha_j.
	\end{equation*}
For any (possible infinity) orthonormal system $\left\{\{\hat{F_j}(\mu,\nu,\alpha)\}_{(\mu,\nu)\in\mathbb{N}^n\times\mathbb{Z}}\right\}_{j\in{J}}$ in $\ell^2(S)$ and any sequence $\{n_j\}_{j\in{J}}$ in $\mathbb{C}$, we have
\begin{equation}\label{Restriction-general}
\bigg\|\sum_{j\in J}n_j|\mathcal{E}_S\{\{\hat{F_j}(\mu,\nu,\alpha)\}|^2\bigg\|_{L^q(\mathbb{T}, L^p(\mathbb{R}_+^n,dw_\alpha))}\leq C\bigg(\sum_{j}|n_j|^{\frac{2p}{p+1}}\bigg)^{\frac{p+1}{2p}},
\end{equation}		
where $C>0$ is independent of $\left\{\{\hat{F_j}(\mu,\nu,\alpha)\}_{(\mu,\nu)\in\mathbb{N}^n\times\mathbb{Z}}\right\}_{j\in{J}}$ and $\{n_j\}_{j\in{J}}$.
\end{theorem}

 {\bf Proof Theorem \ref{main-restriction}.} When $(p,q)=(1, \infty)$, it follows immediately from the fact that the operator $e^{-itL_\alpha}$ is unitary 
 \begin{equation*}
 \bigg\|\sum_{j\in J}n_j|e^{-itL_\alpha}f_j|^2\bigg\|_{L^\infty(\mathbb{T}, L^1(\mathbb{R}_+^n,dw_\alpha))}\leq \sum_{j}|n_j|,
 \end{equation*}		
for any (possible infinity) system $\{f_j\}_{j\in J}$ of orthonormal functions in $L^2(\mathbb{R}_+^n,dw_\alpha)$ and any coefficients $\{n_j\}_{j\in J}$ in $\mathbb{C}$. It is equivalent to 
 \begin{equation*}
 \bigg\|\sum_{j\in J}n_j|\mathcal{E}_S\{\{\hat{F_j}(\mu,\nu,\alpha)\}|^2\bigg\|_{L^\infty(\mathbb{T}, L^1(\mathbb{R}_+^n,dw_\alpha))}\leq \sum_{j}|n_j|,
 \end{equation*}		
 for any (possible infinity) orthonormal system $\left\{\{\hat{F_j}(\mu,\nu,\alpha)\}_{(\mu,\nu)\in\mathbb{N}^n\times\mathbb{Z}}\right\}_{j\in{J}}$ in $\ell^2(S)$ and any sequence $\{n_j\}_{j\in{J}}$ in $\mathbb{C}$.
 
 By Theorem \ref{Restriction-single}, Lemma \ref{duality-principle} and Lemma \ref{general-schatten-bound}, we prove that the inequality \eqref{Restriction-general} holds for $p,q\geq 1$ satisfying
  \begin{equation*}
  1+\frac{1}{n+\sum_{j=1}^{n}\alpha_j}\leq p<\frac{2(n+\sum_{j=1}^{n}\alpha_j)+1}{2(n+\sum_{j=1}^{n}\alpha_j)-1}\quad\text{and}\quad \frac{1}{q}+\dfrac{n+\sum_{j=1}^{n}\alpha_j}{p}=n+\sum_{j=1}^{n}\alpha_j.
  \end{equation*} 
 Combined with $(p,q)=(1,\infty)$ and the diagonal case $(p,q)=(1+\frac{1}{n+\sum_{j=1}^{n}\alpha_j},1+\frac{1}{n+\sum_{j=1}^{n}\alpha_j})$ in Theorem \ref{diagonal}, by interpolation, we obtain our desired estimates.\qed
 \\[11pt]
 
 At the same time, we also obtain the Strichartz inequalities for the system of orthonormal functions.
 \begin{theorem}(Strichartz inequalities for orthonormal functions-general case)\label{main-strichartz}  If $p,q,n\geq1$ such that 
 	\begin{equation*}
 	1\leq p < \dfrac{2(n+\sum_{j=1}^{n}\alpha_j)+1}{2(n+\sum_{j=1}^{n}\alpha_j)-1}\ \ \ \  and \ \ \ \ \dfrac{1}{q}+\dfrac{n+\sum_{j=1}^{n}\alpha_j}{p}=n+\sum_{j=1}^{n}\alpha_j
 	\end{equation*}
 	Then for any (possible infinity) system $\{f_j\}_{j\in J}$ of orthonormal functions in $L^2(\mathbb{R}_+^n,dw_\alpha)$ and any coefficients $\{n_j\}_{j\in J}$ in $\mathbb{C}$, we have
 	\begin{equation}\label{Strichartz-general}
 	\bigg\|\sum_{j\in J}n_j|e^{-itL_\alpha}f_j|^2\bigg\|_{L^q(\mathbb{T}, L^p(\mathbb{R}_+^n,dw_\alpha))}\leq C\bigg(\sum_{j}|n_j|^{\frac{2p}{p+1}}\bigg)^{\frac{p+1}{2p}},
 	\end{equation}		
 	where $C>0$ is independent of $\{f_j\}_{j\in J}$ and $\{n_j\}_{j\in J}$.
 \end{theorem}

The inequality \eqref{Strichartz-general} can be also rewritten in a convenient form in terms of the operator
\begin{equation*}
\gamma_0:=\sum_{j=0}n_j|f_j\rangle_\alpha\langle f_j|
\end{equation*}
which acts on $L^2(\mathbb{R}_+^n,dw_\alpha)$. Here we use the Dirac's notation $|g \rangle_\alpha\langle h|$ for the rank-one operator $f\longmapsto \langle h,f\rangle_\alpha g$. For such $\gamma_0$, define
\begin{equation*}
\gamma(t):=e^{-itL_\alpha}\gamma_0e^{itL_\alpha}=\sum_{j=0}n_j|e^{-itL_\alpha}f_j\rangle_\alpha \langle e^{-itL_\alpha}f_j|.
\end{equation*}
The density of $\gamma(t)$ is given by
\begin{equation*}
\rho_{\gamma(t)}:=\sum_{j=0}n_j|e^{-itL_\alpha}f_j|^2.
\end{equation*}
The inequality \eqref{Strichartz-general} is equivalent to 
\begin{equation*}
\|\rho_\gamma\|_{L^q(\mathbb{T}, L^p(\mathbb{R}_+^n,dw_\alpha))}\leq C\|\gamma_0\|_{\mathcal{G}^{\frac{2p}{p+1}}(L^2(\mathbb{R}_+^n, dw_\alpha))}.
\end{equation*}		

Finally, as an application of Theorem \ref{main-strichartz}, we can obtain the global well-posedness for the Laguerre-Hartree equation in Schatten space using the same method as the proof of Theorem $14$ in \cite{FS}.
\begin{theorem}
Let $w\in L^p(\mathbb{R}_+^n,dw_\alpha)$. Under the same hypotheses as Theorem \ref{main-strichartz}, for any $\gamma_0\in\mathcal{G}^\frac{2p}{p+1}(L^2(\mathbb{R}_+^n, dw_\alpha)$, there exists a unique $\gamma\in C_t^0(\mathbb{T}, \mathcal{G}^\frac{2p}{p+1}(L^2(\mathbb{R}_+^n, dw_\alpha)))$ satisfying $\rho_{\gamma}\in L^q(\mathbb{T},L^p(\mathbb{R}_+^n, dw_\alpha))$ and
\begin{equation*}
\left\{
\begin{array}{rl}
i\partial_t\gamma&=[L_\alpha+w*\rho_\gamma,\gamma]\\
\gamma(0)&=\gamma_0.
\end{array}
\right.
\end{equation*}
\end{theorem}

	{\footnotesize The author declares to have no competing interests.}\\[4mm]
	\noindent\bf{\footnotesize Acknowledgements}\quad\rm
	{\footnotesize The work is supported by the National Natural
		Science Foundation of China (Grant No. 11371036), the Natural Science Foundation of Shaanxi Province (Grant No. 2020JQ-112) and the Fundamental Research Funds for the Central Universities (Grant No. 3102015ZY068).}\\[4mm]
	
	\noindent{\bbb{References}}
	
\end{document}